\documentstyle[12pt]{article}

\begin{titlepage}

\title{Eisenstein series and Scattering matrices}

\author{Ulrich Bunke\thanks{Humboldt-Universit\"at zu Berlin, Institut f\"ur
Reine Mathematik (SFB288), Ziegelstr. 13a, Berlin 10099.
E-mail:ubunke@mathematik.hu-berlin.de
} and
Martin
Olbrich\thanks{Humboldt-Universit\"at zu Berlin, Institut f\"ur Reine
Mathematik (SFB288), Ziegelstr. 13a, Berlin 10099.
E-mail:olbrich@mathematik.hu-berlin.de  }
}

\end{titlepage}

\newcommand{\kaaa}{{\bf k}}

\newcommand{\paaa}{{\bf p}}

\newcommand{\Z}{{\bf Z}}

\newcommand{\C}{{\bf C}}

\newcommand{\gaaa}{{\bf g}}

\newcommand{\aaaa}{{\bf a}}

\newcommand{\naaa}{{\bf n}}

\newcommand{\Hom}{{\mbox{Hom}}}

\newcommand{\End}{{\mbox{ End}}}

\newcommand{\Ree}{{\rm Re }}

\newcommand{\ee}{{\rm e}}

\newcommand{\tr}{{\mbox{tr}}}

\newcommand{\ad}{{\mbox{ad}}}

\newcommand{\id}{{\mbox{id}}}

\newcommand{\supp}{{\mbox{supp}}}

\newcommand{\aca}{{\aaaa_\C^\ast}}

\newtheorem{prop}{Proposition}[section]

\newtheorem{lem}[prop]{Lemma}

\newtheorem{ddd}[prop]{Definition}

\newtheorem{theorem}[prop]{Theorem}

\newtheorem{ass}[prop]{Assumption}

\newtheorem{prob}[prop]{Problem}

\begin{document}

\maketitle

\tableofcontents

\section{The extension problem}\label{sec1}

Let $G$ be a real semisimple Lie group, $G=KAN$ be an Iwasawa decomposition
of $G$, $M:=Z_K(A)$ be the centralizer of $A$ in $K$ and $P=MAN$
be a minimal parabolic subgroup.
Then we have the symmetric space $X=G/K$ and its boundary $\partial
X=G/P=K/M$.
Let $g=\kappa(g)a(g)n(g)$, $\kappa(g)\in K$, $a(g)\in A$, $n(g)\in N$ be
defined with respect to the given Iwasawa decomposition.
Let $\gaaa=\kaaa\oplus\aaaa\oplus\naaa$ be the Iwasawa decomposition of the Lie
algebra $\gaaa$.
If $\lambda\in \aca$, then we set $a^\lambda:=\ee^{\langle
\lambda,\log(a)\rangle}\in\C$.
Corresponding to $\naaa$ there is a positive Weyl chamber $\aaaa^+$.
Let $W=W(\gaaa,\aaaa)$ be the Weyl group generated by the reflections
at walls   of $\aaaa^+$.

We consider a discrete subgroup $\Gamma\subset G$.
\begin{ass}\label{asss}
We assume that there is a $\Gamma$-invariant partition $\partial X =\Omega\cup
\Lambda$ such that
$\Gamma$ acts freely and co-compactly on the open non-empty set $\Omega$.
\end{ass}
We call $\Lambda$ the limit set of $\Gamma$ though this is a slight abuse of
the notion.
For $\lambda,\mu\in\aaaa^*$ we say that $\lambda>\mu$, iff
$\lambda-\mu\in\aaaa^\ast_+$.
Define $\rho\in \aaaa_+^*$ as usual by
$\rho(H):=\frac{1}{2}\tr(\ad(H)_{|\naaa})$, $H\in\aaaa$.
\begin{ddd}
An exponent $\delta_\Gamma\in\aaaa^*$ of $\Gamma$ is a minimal element such
that
$$\sum_{g\in\Gamma} a(gkM)^{-\delta_\Gamma-\lambda-\rho}<\infty$$
converges
uniformly for the parameters $kM$ and $\lambda$ varying in compact
subsets of $\Omega$ and $\aaaa_+^*$, respectively.
\end{ddd}
$\Gamma\subset G$ is discrete and has
a fundamental domain. Since $G$ has at most exponential volume growth
we have clearly
$\delta_\Gamma<\infty$.
Note that we
admit $\delta_\Gamma=-\infty$ and $\delta_\Gamma$ need not to be unique.

Let $\sigma$ be a finite dimensional unitary representation of $M$ on
$V_\sigma$
and $\lambda\in \aca$. Then we can form the representation $\sigma_\lambda$
of $P$ on $V_{\sigma_\lambda}:=V_\sigma$ by
$$\sigma_\lambda(man)=\sigma(m)a^{\rho-\lambda}\ .$$
Let $V(\sigma_\lambda)=G\times_P V_{\sigma_\lambda}$ be the associated
homogeneous bundle.
By $B$ we denote the compact quotient $\Gamma\backslash\Omega$ and let
$V_B(\sigma_\lambda):=\Gamma\backslash V(\sigma_\lambda)$.

Let $\theta$ be the Cartan involution associated to $K$ and set
$\bar{\naaa}=\theta\naaa$.
Then we have the decomposition $\gaaa=\bar{\naaa}\oplus\paaa$, and we can
identify
the tangent bundle of $\partial X$ with $G\times_P \bar{\naaa}$.
Note that $(\Lambda^{max}\bar{\naaa})^{-1/2}=V_{1_0}$,
where $1$ denotes the trivial representation of $M$.
In particular the bundle of half-densities of $\partial X$ is
$V(1_0)$.
Hence  we can define natural parings of sections of
$V(\sigma_\lambda)$ with sections of $V(\tilde{\sigma}_{-\lambda})$,
where $\tilde{\sigma}$ is the dual representation to $\sigma$.

Besides the different globalizations of the principal
series representation of $G$
$$H^{\sigma,\lambda}_*=C^*(\partial X,V(\sigma_\lambda)),\quad
*=\omega,\infty,-\infty,-\omega$$
we consider the following $\Gamma$-modules :
\begin{itemize}
\item $H^{\sigma,\lambda}_{-\omega,\Lambda}:=\{f\in C^{-\omega}(\partial X,
V(\sigma_\lambda)),\supp(f)\subset \Lambda\}$
\item $H^{\sigma,\lambda}_{-\infty}(\Omega):=C^{-\infty}(\Omega,
V(\sigma_\lambda)) $
\item $H^{\sigma,\lambda}_{-\omega,[\Lambda]}:= \{f\in C^{-\omega}(\partial X,
V(\sigma_\lambda)), res(f)\in H^{\sigma,\lambda}_{-\infty}(\Omega)\}$,
\end{itemize}
where $res$ is the restriction to $\Omega$, $-\infty$ stands for distributions
and $-\omega$ for hyperfunctions.
These spaces come with natural topologies.
In fact the first is a Frechet space with the topology induced
from its embedding in the Frechet space of all hyperfunctions.
The second is the topological dual of the space of
smooth functions on $\Omega$ with compact support.
Since the hyperfunctions form a flabby sheaf we have the exact sequence
$$0\rightarrow H^{\sigma,\lambda}_{-\omega,\Lambda} \rightarrow
H^{\sigma,\lambda}_{-\omega,[\Lambda]}\stackrel{res}{\rightarrow}
H^{\sigma,\lambda}_{-\infty}(\Omega) \rightarrow 0$$
inducing a natural topology on the middle space.
\begin{prob}[Extension problem]
Let $\phi\in{}^\Gamma H^{\sigma,\lambda}_{-\infty}(\Omega)$.
Define an invariant extension $ext(\phi)\in{}^\Gamma
H^{\sigma,\lambda}_{-\omega,[\Lambda]}$
as natural as possible such that $\phi=res\circ ext(\phi)$.
\end{prob}
In fact we will define a natural right inverse of $res$ for $\Ree(\lambda)>>0$,
$\lambda\in\aca$.
The idea is then to continue this extension meromorphically.

Let us assume that $\lambda >\delta_\Gamma$. Then we can define a pairing  of
$ {}^\Gamma H^{\sigma,\lambda}_{-\infty}(\Omega)$ with
$H^{\tilde{\sigma},-\lambda}_\infty$ as follows.
Let $f\in H^{\tilde{\sigma},-\lambda}_\infty$. Consider
$f$ as a function on $K$ with values in $V_\sigma$ being right $M$-invariant.
Then
 $(\pi^{\tilde{\sigma},-\lambda}(g)f)(k)
=a^{-\lambda-\rho}(g^{-1}k)f(\kappa(g^{-1}k))$.
By assumption the average
$$\bar{f}:=\sum_{g\in\Gamma} \pi^{\tilde{\sigma},-\lambda}(g) res(f)$$
converges
 in $C^\infty(\Omega,V(\tilde{\sigma}_{-\lambda}))$.
Of course $\bar{f}\in C^\infty(B,V_B(\tilde{\sigma}_{-\lambda}))$
and the map $f\mapsto \bar{f}$ is continuous.
Let $\phi\in C^{-\infty}(B,V_B(\sigma_\lambda))$.
Then $H^{\tilde{\sigma},-\lambda}_\infty\ni f\mapsto
\langle\phi,\bar{f}\rangle$
is  continuous. Thus we can define a distribution
 $ext(\phi)\in H^{\sigma,\lambda}_{-\infty}$
by \begin{equation}\label{exte}\langle ext(\phi),f\rangle :=\langle
\phi,\bar{f}\rangle\ .\end{equation}
It is easy to see that $C^{-\infty}(B,V_B(\sigma_\lambda))\ni \phi\mapsto
ext(\phi)\in H^{\sigma,\lambda}_{-\infty}$ is continuous.

We now discuss the parametrized version.
In order to speak of holomorphic families $\phi_\lambda\in
C^{-\infty}(B,V_B(\sigma_\lambda))$
we must identify the bundles $V_B(\sigma_\lambda)$ for different $\lambda$.
This can be done as follows.
Fix a basis $\mu_1,\dots \mu_r$ of $\aaaa^*$ such that $\mu_i> \delta_\Gamma$.
Let $e_i\in H^{1,\mu_i}_\infty$ be the element given by the constant function
$1$
on $K$. Then $\bar{e_i}$ exists in $C^\infty(B,V_B(1_{\mu_i}))$.
Writing $\lambda=\sum_{i=1}^r\lambda_i\mu_i$ we have
an isomorphism $V_B(\sigma_\lambda)=V_B(\sigma_0)\otimes \prod_{i=1}^r
V_B(1_{\mu_i})^{\lambda_i}$.
Hence we can write $\phi_\lambda=\phi(\lambda)\prod_{i=1}^r
(\bar{e}_i)^{\lambda_i}$
where $\phi(\lambda)\in C^{-\infty}(B,V_B(\sigma_0))$ is now a family
of distributions in a constant bundle.
We define the family $\phi_\lambda$ to be holomorphic, if the family
$\phi(\lambda)$
is a holomorphic family of distributions.
It is easy to check that this notion of holomorphicity does not depend
on the choices.
In a similar way we define holomorphic families $\psi_\lambda\in
H^{\sigma,\lambda}_{-\infty}$.
We identify the bundles $V(\sigma_\lambda)$ for different $\lambda$
using the sections $e_i$.
When we speak of holomorphic families of maps between section spaces
of bundles parametrized by $\lambda$ we will always assume  identifications of
the bundles as above.

It is  not complicated to see that if $\phi_\lambda$ is holomorphic,
then for $\lambda>\delta_\Gamma$ the extension $ext(\phi_\lambda)$ is
holomorphic, too.
Thus we have shown the following
\begin{lem}
For $\lambda>\delta_\Gamma$ the extension problem has a canonical continuous
solution.
Moreover, the extension of a holomorphic family is again a holomorphic family.
\end{lem}

\section{The case of negative $\delta_\Gamma$}

In this section we assume that $\delta_\Gamma$ is negative, i.e.,
$\delta_\Gamma<0$.
We find a meromorphic continuation of the extension defined
in Section \ref{sec1}. The main tool are the scattering matrices
which are closely related to the Knapp-Stein intertwining operators.

Any $w\in W$ can be represented by an element $m_w\in N(M)$, the
normalizer of $M$. If $\sigma$ is a representation of $M$
we define $\sigma^w$ by $\sigma^w(m):=\sigma(m_w^{-1}mm_w)$.
Depending on the choices of $m_w$ there are $G$-invariant
operators
$$\hat{J}_{w,\sigma,\lambda}:H^{\sigma,\lambda}_{*}\rightarrow
H^{\sigma^w,\lambda^w}_{*}$$
which form meromorphic families of pseudodifferential operators
\cite{knappstein71}.
Let $\bar{N}:=exp(\bar{\naaa})$.
Consider $f\in H^{\sigma,\lambda}_{\infty}$ as a right $P$-invariant function
on $G$ with values
in $V_{\sigma_\lambda}$. For $\lambda<0$ the intertwining operator is defined
by the
convergent integral
$$(\hat{J}_{w,\sigma,\lambda}f)(g):=\int_{\bar{N}\cap w^{-1}Nw} f(gm_w\bar{n})
d\bar{n}\ .$$
For other parameters it is obtained by meromorphic continuation.
It is known that $\hat{J}_{w,\sigma,\lambda}$ is compatible
with the choices $*\in\{\omega,\infty,-\infty,-\omega\}$.
Let $c_{w,\sigma}(\lambda)$ be the value of
$\hat{J}_{w^{-1},\sigma^w,\lambda^w}$
on the minimal $K$-type (see \cite{knapp86}, Ch. XV for all that) of the
principal series representation $H^{\sigma,\lambda}_{*}$.
Then $c_{w,\sigma}(\lambda)$ is a meromorphic function on $\aca$ and we define
the normalized intertwining operators by
$$J_{w,\sigma,\lambda}:=
c_{w^{-1},\sigma^w}^{-1}(\lambda^w)\hat{J}_{w,\sigma,\lambda}\ .$$
They satisfy the functional  equations
$$J_{w_1w_2,\sigma,\lambda}=J_{w_1,\sigma^{w_2},\lambda^{w_2}}\circ
J_{w_2,\sigma,\lambda},\quad w_1,w_2\in W\ .$$
A special case is
\begin{equation}\label{spex}J_{w^{-1},\sigma^w,\lambda^w}\circ
J_{w,\sigma,\lambda}=\id,\quad w\in W\ .\end{equation}

In order to simplify the notation we assume that $\sigma$ is Weyl
invariant, i.e. $\sigma^w=\sigma$. This can always be achieved by taking
the direct sum of all Weyl translates of $\sigma$.

The scattering matrices are operators
$$S_{w,\sigma,\lambda}:C^{-\infty}(B,V_B(\sigma_\lambda))\rightarrow
C^{-\infty}(B,V_B(\sigma_{\lambda^w}))\ .$$
\begin{ddd}\label{scatdef}
For $\lambda>\delta_\Gamma$ we define
$$S_{w,\sigma,\lambda}(\phi)=res\circ J_{w,\sigma,\lambda}\circ ext(\phi)\ .$$
\end{ddd}
Then the $S_{w,\sigma,\lambda}$, $w\in W$, are  meromorphic families of
pseudodifferential operators.
Locally they coincide with $J_{w,\sigma,\lambda}$ up
to  smoothing operators.
Since $\delta_\Gamma<0$ the scattering matrices are defined
on a neighbourhood  $C$ of $\{\lambda\in \aca | \Ree (\lambda)\ge 0\}$.
They satisfy the functional  equations
\begin{equation}\label{funeq}
S_{w_1w_2,\sigma,\lambda}=S_{w_1,\sigma,\lambda^{w_2}}\circ
S_{w_2,\sigma,\lambda},\quad w_1,w_2\in W\ ,
\end{equation}
when all terms are defined.

We now provide a simultaneous meromorphic continuation of $ext$ and the
scattering matrices
to all of $\aca$.
For $\lambda\in C$ we define
\begin{equation}\label{open}S_{w^{-1},\sigma,\lambda^w}:=
S^{-1}_{w,\sigma,\lambda}\ .\end{equation}
We claim that this formula defines a meromorphic family of operators.

To see the claim
we invoke the Fredholm theory for Frechet spaces \cite{grothendieck56}.
Since $\Gamma$ acts properly on $\Omega$ there is a $\Gamma$-invariant
function
$\chi\in C^\infty(\Omega\times\Omega)$ being identically one
on a small neighbourhood of the diagonal and zero outside  of a somewhat
larger neighbourhood. We assume the latter neighbourhood to be so small that
$\chi(x,gx)=0$ for all $x\in\Omega$, $1\not=g\in \Gamma$.
Let $\tilde{J}_{w,\sigma,\lambda}$ be the operator obtained by cutting
off the distributional kernel of $J_{w,\sigma,\lambda}$ with $\chi$.
Then we can consider $\tilde{J}_{w,\sigma,\lambda}$ as a meromorphic
family of operators on $B$.
Since the singularities of the (unnormalized) intertwining operators
are local operators
we conclude  from (\ref{spex}) that $\tilde{J}_{w^{-1},\sigma,\lambda^w}\circ
S_{w,\sigma,\lambda}=\id+ R_1(\lambda)$,
$S_{w,\sigma,\lambda}\circ \tilde{J}_{w^{-1},\sigma,\lambda^w}=
\id+R_2(\lambda)$,
where $R_1,R_2$ are holomorphic families of smoothing
operators. Thus
\begin{equation}\label{finit}S_{w,\sigma,\lambda}^{-1}=(\id+
R_1(\lambda))^{-1}\circ \tilde{J}_{w^{-1},\sigma,\lambda^w}\end{equation}
is meromorphic if $(\id+ R_1(\lambda))^{-1}$ is so.
The meromorphicity of the latter family follows from the Fredholm theory for
Frechet spaces  \cite{grothendieck56}
if $(\id+ R_1(\lambda))$ is injective for some $\lambda\in C$.
But there is an open subset in $C\cap C^{w^{-1}}$ where (\ref{open}) holds.

We call $\lambda\in \aca$ bad if
$$2 \frac{\langle \lambda,\alpha\rangle
}{\langle\alpha,\alpha\rangle}\in\Z\quad  \mbox{for some root $\alpha$ of
$(\gaaa,\aaaa)$}\ ,$$
 and good otherwise,
where $\langle.,.\rangle$ is some invariant $\C$-bilinear scalar product on
$\aca$.
It is known that the normalized intertwining operators are holomorphic
at  good $\lambda$. From (\ref{finit}) we obtain that at good points
$\lambda_0\in\aca$
the scattering matrix $S_{w,\sigma,\lambda}^{-1}$
has at most finite dimensional singularities
in the following sense. There exists a family of operators $P(\lambda)$,
holomorphic near $\lambda_0$ and having at most finite dimensional kernels
for $\lambda$ near $\lambda_0$ such that the composition $P(\lambda)\circ
S_{w,\sigma,\lambda}^{-1}$ is holomorphic. In fact take $P(\lambda):= \id+
R_1(\lambda) $.
(In this way we will understand the finite dimensionality of the singularities
of $ext$ and the Eisenstein series later on.)

This way we have continued $S_{w^{-1},\sigma,\lambda}$ to $C^w$ for all $w\in
W$.
Next we define the extension of  $\phi\in C^{-\infty}(B,V_B(\sigma_\lambda))$,
$\lambda\in C^w$, by

$$ext(\phi):=J_{w,\sigma,\lambda^{w^{-1}}}\circ ext\circ
S_{w^{-1},\sigma,\lambda}(\phi)\ .$$
Since $\bigcup_{w\in W}C^w=\aca$ the extension  $ext$  becomes a meromorphic
family of continuous maps $ext: C^{-\infty}(B,V(\sigma_\lambda))\rightarrow
{}^\Gamma H^{\sigma,\lambda}_{-\infty}$ on all of $\aca$.
In order to discuss the singularities of $ext$ we consider the representation
$$S_{w^{-1},\sigma,\lambda}=\tilde{J}_{w^{-1},\sigma,\lambda}\circ (\id+
R_2(\lambda))^{-1}\ .$$
 Then $$ext=J_{w,\sigma,\lambda^{w^{-1}}}\circ ext\circ
\tilde{J}_{w^{-1},\sigma,\lambda}\circ (\id+  R_2(\lambda))^{-1}\ .$$
The local singularities of the intertwining operators at bad points cancel,
and thus
$J_{w,\sigma,\lambda^{w^{-1}}} \linebreak[4]\circ ext\circ
\tilde{J}_{w^{-1},\sigma,\lambda}$
is holomorphic on all of $C^w$.
 $(\id+  R_2(\lambda))^{-1}$ has finite dimensional singularities, and so has
$ext$.

Now we can define
all scattering matrices on $\aca$ by \ref{scatdef}
obtaining meromorphic families of pseudodifferential operators.
At good points they have at most finite dimensional singularities.
We have  shown
\begin{prop}
Assume $\delta_\Gamma<0$. Then there is a meromorphic continuation of
$ext:C^{-\infty}(B,V(\sigma_\lambda))\rightarrow {}^\Gamma
H^{\sigma,\lambda}_{-\infty}$
to all of $\aca$. The singularities of $ext$ are finite dimensional.
Moreover there are
meromorphic continuations of
the scattering matrices satisfying the functional equations (\ref{funeq}). At
good points
the singularities are  at most finite  dimensional.
\end{prop}

\section{The general case}

In this subsection we assume that $G$ belongs to a series of equal real rank.
This means that there is a sequence $\dots \subset G^n\subset G^{n+1}\subset
\dots$
of real semisimple Lie groups inducing embeddings of the corresponding
Iwasawa constituents $K^n\subset K^{n+1}$,
$N^n  \stackrel{\scriptstyle  \subset}{ \scriptstyle \not=}  N^{n+1}$, $
M^n\subset M^{n+1}$
such that $A=A^n = A^{n+1}$.
Then we have totally geodesic embeddings of the symmetric spaces
$X^n\subset X^{n+1}$
inducing embeddings of their boundaries
$\partial X^n\subset \partial X^{n+1}$.
If $\Gamma\subset G^n$ satisfies \ref{asss} then it keeps satisfying
\ref{asss}
when viewed  as a subgroup of $G^{n+1}$.
We obtain the embedding
$\Omega^n\subset  \Omega^{n+1}$ inducing
$B^n\subset  B^{n+1}$
while the limit set $\Lambda^n$ is identified with $\Lambda^{n+1}$.
Let $\rho^n(H)=\frac{1}{2}\tr(\ad(H)_{|\naaa^n})$, $H\in\aaaa$.

The exponent of $\Gamma$ now depends on $n$ and is denoted by
$\delta_\Gamma^n$.
We have the relation $\delta_\Gamma^{n+1}=\delta_\Gamma^n-\rho^{n+1}+\rho^n$.
Hence $\delta_\Gamma^{n}\to -\infty$ as $n\to\infty$. In  particular if
we embed $\Gamma\subset G^n\subset G^{n+m}$ for large enough $m$
we can obtain $\delta_\Gamma^{n+m}<0$ and solve the extension problem as well
as
the continuation of the scattering matrices.

Let $\Gamma\subset G^n$ satisfy \ref{asss}.
The aim of the present section is to show how the solution $ext^{n+1}$ of
the extension problem on $\partial X^{n+1}$ leads to the solution $ext^n$ of
the extension
problem on $\partial X^n$.

Let $\sigma$ be a Weyl invariant representation of $M^{n+1}$. Then it
restricts to a Weyl invariant representation of $M^n$.
The representation $\sigma_\lambda$ of $P^{n+1}$ restricts to the
representation
$\sigma_{\lambda+\rho^n-\rho^{n+1}}$ of $P^n$.
This induces an isomorphism of bundles
 $$V_{B^{n+1}}(\sigma_\lambda)_{|B^n}=
V_{B^n}(\sigma_{\lambda+\rho^n-\rho^{n+1}})\ .$$
We obtain a push forward of distributions
$$i_*:C^{-\infty}(B^n,V_{B^n}(\sigma_\lambda))\rightarrow
C^{-\infty}(B^{n+1},V_{B^{n+1}}(\sigma_{\lambda+\rho^n
-\rho^{n+1}}))\ .$$
Assume that we have a solution $ext^{n+1}$ of the extension problem.
For\linebreak[4] $\phi\in C^{-\infty}(B^n,V_{B^n}(\sigma_\lambda))$ the push
forward
$i_\ast(\phi)$ has support in $\partial X^n$. Thus $ext(\phi)$
has support in $\partial X^n$, too.
Let $f\in C^\infty(\partial
X^{n+1},V(\tilde{\sigma}_{-\lambda-\rho^n+\rho^{n+1}}))$
such that $f_{|\partial X^n}=0$.
We claim that $\langle ext^{n+1}\circ i_*(\phi),f\rangle=0$.
To see the claim embed $\phi$ into a holomorphic family $\phi_\lambda\in
C^{-\infty}(B^n,V_{B^n}(\sigma_\lambda))$.
Then $\langle ext^{n+1}\circ i_*(\phi_\lambda),f\rangle =0$ for $\lambda>>0$.
Now the claim follows by meromorphic continuation.
Thus we can define a pull back
$ext^n(\phi):=i^*\circ ext^{n+1}\circ i_*(\phi)$ as follows.
For $f\in C^\infty(\partial X^{n},V(\tilde{\sigma}_{-\lambda}))$ let
$\tilde{f}\in C^\infty(\partial
X^{n+1},V(\tilde{\sigma}_{-\lambda-\rho^n+\rho^{n+1}}))$
be an arbitrary extension. We set
$$\langle ext^n(\phi),f\rangle := \langle ext^{n+1}\circ
i_*(\phi_\lambda),\tilde{f}\rangle\ .$$
Then $ext^n$ is well defined.
The meromorphic continuation of $ext$ immediately implies the meromorphic
continuation
of the scattering matrices by \ref{scatdef}.

For any given finite dimensional representation of $M^n$ we can find
a Weyl invariant representation $\sigma$ of $M^{n+1}$ such that $\sigma_{|M^n}$
contains the former as a subrepresentation.
\begin{theorem}
Let $\Gamma\subset G$ satisfy assumption \ref{asss} and assume that $G$ belongs
 to a series.
Then the extension problem has a unique continuous meromorphic solution
$$ext:C^{-\infty}(B,V_{B}(\sigma_\lambda)) \rightarrow {}^\Gamma
H^{\sigma,\lambda}_{-\infty}$$
having finite dimensional singularities.
Moreover,
there are meromorphic continuations of the scattering matrices to
all of $\aca$. They form meromorphic families of pseudodifferential operators
locally coinciding with the intertwining operators up to smoothing operators
and satisfying
the  functional  equations (\ref{funeq}). At good points the singularities
of the  scattering matrices are at most finite dimensional.
\end{theorem}
Note that the singular part of $ext(\phi_\lambda)$ at $\lambda=\lambda_0$
provides $\Gamma$-invariant distributions in $V(\sigma_{\lambda_0})$ with
support in $\Lambda$.
The singularities of the scattering matrices at good points
are often called resonances while at the bad points in general
there is a mixing of the topological singularity localized on the diagonal
and resonances.

\section{Eisenstein series}

In this section we define the Eisenstein series and discuss
its functional equations.
We assume that the extension problem for $\partial  X$ has a meromorphic
solution $ext$.

Let $\sigma$ be a Weyl-invariant finite dimensional representation of $M$ and
$\gamma$ be a finite dimensional representation of $K$.
Let $V(\gamma):=G\times_KV_\gamma$ be the associated homogeneous bundle over
$X$.
The Eisenstein series $E(\lambda,\phi,T)$
depends on the data $\lambda\in \aca$, $\phi\in
C^{-\infty}(B,V_B(\sigma_\lambda))$ and
$T\in\Hom_M(V_\sigma,V_\gamma)$.
There is a Poisson transform
$$P^T_\lambda: H^{\sigma,\lambda}_{-\omega}\rightarrow C^\infty(X,V(\gamma))
$$
 defined by
$$(P^T_\lambda f)(g):=\int_K
a(g^{-1}k)^{-(\lambda+\rho)}\gamma(\kappa(g^{-1}k)) T f(k) dk\ .$$
The range of $P^T_\lambda$ can be characterized by invariant differential
equations
formally written as $D_\lambda P^T_\lambda =0$
\cite{olbrichdiss}, \cite{bunkeolbrich947}.
\begin{ddd}
The Eisenstein series is defined by
$$E(\lambda,\phi,T):= P^T_\lambda \circ ext(\phi)\ .$$
\end{ddd}
Thus $E(\lambda,\phi,T)$ is a meromorphic family of smooth sections of
$V(\gamma)$
satisfying $D_\lambda E(\lambda,\phi,T)\linebreak[4]=0$. Using (\ref{exte}) one
can check that for $\lambda>\delta_\Gamma$
our definition coincides with other definitions of Eisenstein series
in the literature (up to normalizations).
The singularities of the Eisenstein series
are finite dimensional.

We now discuss the functional equations satisfied by the Eisenstein series.
For $\lambda<0$, $w\in W$ define
$$c_{w,\gamma}(\lambda):=\int_{\bar{N}\cap
w^{-1}Nw}a(\bar{n})^{-(\lambda+\rho)}\gamma(\kappa(\bar{n}))d\bar{n}\in
\End_M(V_\gamma)\ .$$
$c_{w,\gamma}(\lambda)$ extends  meromorphically to all of $\aca$.
We employ the functional equations of the Poisson transform
proved in \cite{olbrichdiss}:
$$c_{w,\sigma}(\lambda) P^T_\lambda \circ J_{w^{-1},\sigma,\lambda^w} =
P^{\gamma(m_w)c_{w,\gamma}(\lambda)T}_{\lambda^w}, \quad \forall w\in W\ .$$
It follows
\begin{prop}
The Eisenstein series satisfies the functional equations
$$
 E(\lambda,c_{w,\sigma}(\lambda)S_{w^{-1},\sigma,\lambda^w}\phi,T)
=E(\lambda^w,\phi,\gamma(m_w)c_{w,\gamma}(\lambda)T),\quad \forall w\in W\ .$$
\end{prop}
In the rank one case observe that for $\lambda_0>0$ the singular part of
$E(\lambda,\phi,T)$ at
$\lambda=\lambda_0$ provides $L^2$-solutions of the system $D_\lambda f=0$ on
$Y=\Gamma\backslash X$.
In fact we obtain the singular part by the Poisson transform
of the singular part of $ext(\phi_\lambda)$ which has support in the limit set
$\Lambda$.\\[1cm]
{\bf Historical remarks}\\
In the case that $X$ is the real hyperbolic space of dimension
$n$ and $\gamma$ is the trivial $K$-type, the meromorphic continuation
of the Eisenstein series was previously obtained by
Patterson \cite{patterson75}, \cite{patterson76}, \cite{patterson761},
\cite{patterson89},
Mandouvalos \cite{mandouvalos86}, \cite{mandouvalos88}, \cite{mandouvalos89},
and Perry \cite{perry87}, \cite{perry89}.
The idea of using the fact that $G$ belongs to a series we learned from
Mandouvalos.
As Patterson remarked it can be traced back to Reshnikov and even Selberg.
That the Fredholm theory can be applied to the meromorphic continuation of the
scattering matrices was
first observed by Patterson \cite{patterson76}.
Note that the techniques invoked in the literature above are {\em much}
more complicated than the methods presented in the present paper. The former
are based on the analysis of the resolvent kernel
of the Laplace operator on $Y$. It is a challenging problem to
do this analysis for bundles even over rank one spaces.
In \cite{bunkeolbrich951} we considered the extension problem
in connection with the $\Gamma$-cohomology of
$H^{\sigma,\lambda}_{-\omega,\Lambda}$
for Fuchsian groups of the second kind.

\bibliographystyle{plain}

\end{document}